\documentclass{article}
\usepackage{amsmath,amstext, verbatim,remreset}
\usepackage[T2A]{fontenc}
\usepackage[cp1251]{inputenc}
\usepackage[english]{babel}
\usepackage{amsfonts,amssymb}
\newtheorem{theorem}{Theorem}
\newtheorem{lemma}{Lemma}
\begin{document}
\begin{center}
{\Large \bf Comparison theorems and some of their applications}
\vskip 5 mm
{\Large \bf V.~F.~Babenko, O.~V.~Kovalenko}
\vskip 5 mm
{Oles Gonchar Dnipropetrovsk National University
 \\{\it E-mail: babenko.vladislav@gmail.com }
{\it E-mail: olegkovalenko90@gmail.com} }
\end{center}
\begin{abstract}
{Analogues of Kolmogorov comparison theorems and some of their applications were established.}
\end{abstract}

\section{Notations. Statement of the problem. Known results.} 
Let $L_\infty({\mathbb R})$ denote the space of all measurable and essentially bounded functions $x\colon {\mathbb R}\to {\mathbb R}$ with the norm
$$\|x\|=\|x\|_{L_\infty({\mathbb R})}={\rm
ess\,sup}\left\{|x(t)|:t\in{\mathbb R}\right\}.$$ For natural $r$ let
$L_\infty^r({\mathbb R})$ denote the space of functions $x\colon {\mathbb R}\to
{\mathbb R}$ such that the derivative $x^{(r-1)}, \; x^{(0)}=x,$ is locally absolute continuous and  $x^{(r)}\in L_\infty({\mathbb R})$. Set
$L_{\infty,\infty}^r({\mathbb R}):=L_\infty^r({\mathbb R})\bigcap L_\infty({\mathbb R})$.

For $r\in {\mathbb N}$ by $\varphi_r(t)$ we will denote the Euler prefect spline of the order $r$ (i.~e. $r$-th periodic integral of the functions ${\rm sgn}\sin t$ with zero mean value on the period).
For $\lambda> 0$ set $\varphi_{\lambda, r}(t):=\lambda^{-r}\varphi_{r}(\lambda t).$

To prove his outstanding inequality (see \cite{Kol1}---\cite{Kol3}) Kolmogorov proved a statement, known as  a comparison theorem.

{\bf Theorem A. }{\it Let $r\in{\mathbb N}$ and a function $x\in L^r_{\infty,\infty}({\mathbb R})$ are given. Let numbers $a\in{\mathbb R}$ and $\lambda>0$ be such, that
$$
\|x^{(k)}\|\leq \|a\varphi_{\lambda,r}^{(k)}\|,\,\, k\in\{0,r\}.
$$ If points $\xi,\eta\in{\mathbb R}$ are such, that $x(\xi)=a\varphi_{\lambda,r}(\eta)$, then $$|x'(\xi)|\leq |a|\cdot|\varphi'_{\lambda,r}(\eta)|.$$ }

Both the Kolmogorov comparison theorem and its proof played important role in exact solutions of many extremal problems in approximation theory (see.~\cite{korn1, korn2, korn3, BKKP}).

The goal of this paper is to prove several analogues of Kolmogorov comparison theorems.

In the next paragraph we will introduce a family of splines, which will play the same role, as Euler perfect splines play in the theorem A, and study some of their properties. In \S~3 we will prove 3 analogues of Kolmogorov comparison theorem for the cases when the norms of a function and its derivatives of orders $r-1$ and $r$ are given; the norms of a function and its derivatives of orders  $r-2$ and $r$ are given; the norms of a function and its derivatives of orders $r-2$, $r-1$ and $r$ are given. In \S~4 we will give several applications of the obtained comparison theorems.

\medskip

\section{Comparison functions and their properties.} Let $a_1, a_2\geq 0$.  Set $T:=a_1+a_2+2$. Define a function $\psi_1(a_1,a_2;t)$ in the following way. On the segment $[0,T]$
set
 $$
           \psi_1(a_1,a_2;t):=\left\{
           \begin{array}{rcl}
            0, &  t\in[0,a_1], \\
              t-a_1,   & t\in[a_1,a_1+1], \\
              1, & t\in [a_1+1,a_1+a_2+1], \\
              2+a_1+a_2-t, & t\in [a_1+a_2+1, T].\\
           \end{array}
           \right.
  $$
  Continue function $\psi_1(a_1,a_2;t)$ to the segment $[T,2\cdot T]$ by the equality
  \begin{equation}\label{2.0}
  \psi_1(a_1,a_2;t) = - \psi_1(a_1,a_2;t-T),\,t\in[T,2\cdot T]
  \end{equation} and then periodically with period $2\cdot T$ to the whole real line.

  Note, that $\psi_1(a_1,a_2;t)\in L_{\infty,\infty}^1({\mathbb R})$.

  For $r\in {\mathbb N}$ denote by  $\psi_r(a_1,a_2;t)$ $(r-1)$-th $2\cdot T$ -- periodic integral of the function $\psi_1(a_1,a_2;t)$ with zero mean on the period (so that, particularly,
  $\psi'_r(a_1,a_2;t)=\psi_{r-1}(a_1,a_2;t)$).
Rodov~\cite{Rodov1946} was the first, who considered the functions $\psi_r(a_1,a_2;t)$.

  We will list several properties of the function $\psi_r(a_1,a_2;t)$, $r\in{\mathbb N}$, which can easily be proved either directly from definition, or similar to the corresponding  properties of Euler splines $\varphi_r$ (see, for example,~\cite[Chapter 5]{korn1},~\cite[Chapter 3] {korn2}). Note, that the function $\psi_r$ is $2\cdot T$-periodic and for all $r\geq 1$ 
   $$\psi_r(a_1,a_2;t) = - \psi_r(a_1,a_2;t-T),\,t\in[T,2\cdot T].$$
  
    Moreover, the function $\psi_2(a_1,a_2;t)$ has exactly two zeroes on the period --- the points $a_1+\frac{a_2}{2}+1$ and $2a_1+\frac{3a_2}{2}+3$.
    Hence the functions $\psi_r(a_1,a_2;t)$ for $r\geq 2$ also have exactly two zeroes on the period: for any $k\in {\mathbb N}$
    \begin{equation}\label{1}
    \psi_{2k+1}(a_1,a_2;0)=\psi_{2k+1}(a_1,a_2;a_1+a_2+2)=0,
    \end{equation}
        \begin{equation}\label{2}
     \psi_{2k}\left(a_1,a_2;a_1+\frac{a_2}{2}+1\right)=\psi_{2k}\left(a_1,a_2;2a_1+\frac{3a_2}{2}+3\right)=0.
    \end{equation}

Note, that for $a_1 = 0$ the equality~\eqref{1} is true for $k=0$ too.

Hence, in turn, we have that for $r\geq 3$ (in the case $a_1=0$ for $r\geq 2$) the function $\psi_r(a_1,a_2;t)$ is strictly monotone between zeroes of its derivative and the plot of the function $\psi_r(a_1,a_2;t)$ is strictly convex at the intervals of constant sign. Moreover, it is easy to see, that the plot of the function $\psi_r(a_1,a_2;t)$ is symmetrical with respect to its zeroes and the lines  $t=t_0$, where $t_0$ -- is the zero of  $\psi'_r(a_1,a_2;t)$.
 At last note, that $\psi_r(0,0;t)=\varphi_{\pi/2, r}(t)$.

For  $r\in{\mathbb N}$, $a_1,a_2\geq 0,\,\lambda>0$ and  $b\in{\mathbb R}$ set $$\Psi_{a_1,a_2,b,\lambda}(t)=\Psi_{r;a_1,a_2,b,\lambda}(t):=$$
$$=b\left(\frac{\lambda}{2a_1+2a_2+4}\right)^{r}\psi_r\left(a_1,a_2;\frac{2a_1+2a_2+4}{\lambda} t\right).$$ Note, that the function $\Psi_{a_1,a_2,b,\lambda}(t)$ is $\lambda$ -- periodic.

 \begin{theorem}\label{th::0}
 Let $r\in{\mathbb N}$ and $x\in L^r_{\infty,\infty}({\mathbb R})$ be given. Then
 \begin{enumerate}
 \item[$a)$]   there exist $a_2\geq 0,\,\lambda>0$ and $b\in{\mathbb R}$ such, that
 $$ \left\|\Psi_{0,a_2,b,\lambda}^{(s)}\right\|=\left\|x^{(s)}\right\|,\,\,s\in\left\{0,r-1,r\right\}.$$
  \item[$b)$]   there exist $a_1\geq 0,\,\lambda>0$ and $b\in{\mathbb R}$ such, that
  $$ \left\|\Psi_{a_1,0,b,\lambda}^{(s)}\right\|=\left\|x^{(s)}\right\|,\,\,s\in\left\{0,r-2,r\right\}.$$
   \item[$c)$]  there exist $a_1,a_2\geq 0,\,\lambda>0$ and $b\in{\mathbb R}$ such, that
    $$ \left\|\Psi_{a_1,a_2,b,\lambda}^{(s)}\right\|=\left\|x^{(s)}\right\|,\,\,s\in\left\{0,r-1,r-2,r\right\}.$$
 \end{enumerate}
\end{theorem}

The truth of the theorem~\ref{th::0}, actually, follows from Kolmogorov comparison theorem.  We will prove the statement $a)$, the rest of the statements can be proved analogously.
 
 It is clear, that $\Psi_{0,0,b,\lambda}(t)=b\left(\frac\lambda 4\right)^r\varphi_{\pi/2}(\frac 4\lambda t)$. Hence the parameters $b$ and $\lambda$ can be chosen in such way, that $ \left\|\Psi_{0,0,b,\lambda}^{(s)}\right\|=\left\|x^{(s)}\right\|$, $s=r-1,r$. Then theorem~A implies that $ \left\|\Psi_{0,0,b,\lambda}\right\|\leq\left\|x\right\|$. When the parameter $a_2$ increases  $ \left\|\Psi_{0,a_2,b,\lambda}\right\|$ continuously increases from $ \left\|\Psi_{0,0,b,\lambda}\right\|$ to $\infty$, and $ \left\|\Psi_{0,a_2,b,\lambda}^{(s)}\right\|,\; s=r-1, r$ remain unchanged. Hence we can choose the parameter $a_2$ such, that $ \left\|\Psi_{0,a_2,b,\lambda}^{(s)}\right\|=\left\|x^{(s)}\right\|$, $s=0,r-1,r$.


\section{Comparison theorems.}

The next theorem contains three analogues of Kolmogorov comparison theorem.
 \begin{theorem}\label{th::1}
 Let $r\in{\mathbb N}$ and $x\in L_{\infty,\infty}^r({\mathbb R})$ be given. Let one of the following condition holds.
\begin{enumerate}
\item[$a)$] The numbers $a_1 = 0$, $a_2\geq 0$, $\lambda>0$ and $b\neq 0$ are such, that
 \begin{equation}\label{4.1}
\left\|x^{(k)}\right\|\le
\left\|\Psi^{(k)}_{a_1,a_2,b,\lambda}\right\|,\,\,k\in\{0,r-1,r\}.
 \end{equation}
\item[$b)$] The numbers $a_1 \geq 0$, $a_2= 0$, $\lambda>0$ and $b\neq 0$ are such, that
 \begin{equation}\label{4.2}
\left\|x^{(k)}\right\|\le
\left\|\Psi^{(k)}_{a_1,a_2,b,\lambda}\right\|,\,\,k\in\{0,r-2,r\}.
 \end{equation}
\item[$c)$] The numbers $a_1 \geq 0$, $a_2\geq 0$, $\lambda>0$ and $b\neq 0$ are such, that
 \begin{equation}\label{4.3}
\left\|x^{(k)}\right\|\le
\left\|\Psi^{(k)}_{a_1,a_2,b,\lambda}\right\|,\,\,k\in\{0,r-1,r-2,r\}.
 \end{equation}
\end{enumerate}

  If points $\tau$ and $\xi$ are such that $x(\tau) = \Psi_{a_1,a_2,b,\lambda}(\xi)$, then
\begin{equation}\label{6'}\left|x'(\tau)\right|\le\left|\Psi'_{a_1,a_2,b,\lambda}(\xi)\right|.
 \end{equation}
 \end{theorem}

{\bf Proof.} For brevity we will write $\Psi (t)$ instead of $\Psi_{a_1,a_2,b,\lambda}(t)$ in the proof of this theorem. Considering, if necessary, the function $-x(t)$ instead of $x(t)$ and function $-\Psi (t)$ instead of 
$\Psi (t)$, we can count that $x'(\tau)>0$ and
 \begin{equation}\label{7}
 \Psi'(\tau)>0.
 \end{equation}
Moreover, considering appropriate shift $\Psi (\cdot +\alpha)$ of the function $\Psi$, we can count that $\tau = \xi$, i.~e.
 \begin{equation}\label{6}
 x(\tau)=\Psi(\tau).
 \end{equation}
Assume, that~\eqref{6} holds, but instead the inequality~(\ref{6'}) (with $\xi =\tau$) the inequality
$$\left|x'(\tau)\right|>\left|\Psi'(\tau)\right| $$ holds. Denote by $(\tau_1,\tau_2)$ the smallest interval which contains  $\tau$ on which the function $\Psi$ is monotone and such that $\Psi'(\tau_1)=\Psi'(\tau_2)=0$.
In virtue of the assumption there exists a number $\delta>0$ such that $x'(t)>\Psi'(t)$ for all $t\in (\tau-\delta,\tau+\delta)$, and hence in virtue of~\eqref{6} $x(\tau+\delta)>\Psi(\tau+\delta)$ and $x(\tau-\delta)<\Psi(\tau-\delta)$.

Choose $\varepsilon>0$ so small, that for a function
$x_\varepsilon(t):=(1-\varepsilon)x(t)$ the following inequalities hold:
$x_\varepsilon(\tau+\delta)>\Psi(\tau+\delta)$ and 
$x_\varepsilon(\tau-\delta)<\Psi(\tau-\delta)$. In virtue of the conditions~\eqref{4.1} -- \eqref{4.3} and condition 
\eqref{7} we have
$$x_\varepsilon(\tau_1)>\Psi(\tau_1),\;\; x_\varepsilon(\tau_2)<\Psi(\tau_2).$$
Hence on the interval $(\tau_1,\tau_2)$ the difference
$\Delta_\varepsilon(t):=x_\varepsilon(t)-\Psi(t)$ has at least 3 sign changes.

It is easy to see, that there exist a sequence of functions $\mu_N\in
C^\infty({\mathbb R})$, $N\in {\mathbb N}$ with the following properties:
\begin{enumerate}
\item $\mu_N(t)=1$ on interval $[\tau_1,\tau_2]$; $\|\mu_N\|=1$;
\item $\mu_N(t)=0$ for all $t$ outside the interval
$[\tau_1-N\cdot\frac{\lambda}2;\tau_1+N\cdot\frac{\lambda}2]$ (where, as before, $T=a_1+a_2+2$);
\item for all
$k=1,2,\dots,r$
$$\max\limits_{j=\overline{1,k}}\left\|\mu_N^{(j)}\right\|<
\varepsilon \|x_\varepsilon^{(k)}\|\left(\sum\limits_{i=1}^k
C_k^i\left\|x_\varepsilon^{(k-i)}\right\|\right)^{-1},$$ if $N$
is enough big.
\end{enumerate}

Below we count that $N$ is chosen enough big, so that the property 3 holds.

Set $$x_N(t):=x_\varepsilon(t)\cdot\mu_N(t),$$
and
$$\Delta_N(t):=\Psi(t)-x_N(t).$$ Then
$$x_N(t)=x_\varepsilon(t),\;if\; t\in[\tau_1,\tau_2],$$
\begin{equation}\label{7'}
\Delta_N(t)=\Psi(t),if\,|t-\tau_1|\geq N\cdot\frac\lambda 2
\end{equation}
and
$$
\|x_N\|\leq\|x_\varepsilon\|=(1-\varepsilon)\|x\|\leq(1-\varepsilon)\|\Psi\|.
$$
Moreover, for $k=1,\dots,r$ $$\left|x_N^{(k)}(t)\right|=
\left|\left[x_\varepsilon(t)\mu_N(t)\right]^{(k)}\right|=
\left|\sum\limits_{i=0}^k C_k^i
x_\varepsilon^{(k-i)}(t)\mu_N^{(i)}(t)\right|\leq
$$
$$
\leq \left\|x_\varepsilon^{(k)}\right\|+\sum\limits_{i=1}^k
C_k^i\left\|x_\varepsilon^{(k-i)}\right\|\left\|\mu_N^{(i)}\right\|.$$ Hence,
in virtue of property 3 of the function $\mu_N$ and the choice of the number $N$, we get
$$
\left\|x_N^{(k)}\right\|<
\left\|x_\varepsilon^{(k)}\right\|+\varepsilon\left\|x^{(k)}_\varepsilon\right\|=(1-\varepsilon)\left\|x^{(k)}\right\|+\varepsilon\left\|x^{(k)}\right\|=\left\|x^{(k)}\right\|.
$$
For $t\in [\tau_1,\tau_2]$ we have
$\Delta_N=\Psi(t)-x_\varepsilon(t),$ and hence the function $\Delta_N(t)$
has at least three sign changes on the interval $[\tau_1,\tau_2]$. At each of the rest monotonicity intervals of the function $\Psi$ the function 
$\Delta_N$ has at least one sign change. Hence on the interval $\left[\tau_1-N\cdot\frac{\lambda}2,\tau_1+N\cdot\frac{\lambda}2\right]$
the function $\Delta_N(t)$ has at least $2N+2$ sign changes. Moreover, in virtue of~\eqref{1},\eqref{2} and \eqref{7'} for all
$i=1,2,\dots,\left[\frac{r-1}{2}\right]$ the following equalities hold
\begin{equation}\label{10}
\Delta_N^{(2i-1)}\left(\tau_1-N\cdot\frac\lambda 2\right)=\Delta_N^{(2i-1)}
\left(\tau_1+N\cdot\frac\lambda 2\right)=0.
\end{equation}

All of the arguments above are true if any of the condition $a)-c)$ hold.
Let now condition $a)$ of the theorem holds.

Applying Rolle's theorem and counting~\eqref{10} we have that the function $\Delta_N^{(r-1)}(t)$ has at least $2N+2$ zeroes on the interval $$\left[\tau_1-N\cdot\frac{\lambda} 2,\tau_1+N\cdot\frac{\lambda} 2\right].$$
Hence on some monotonicity interval
$$[\alpha,\alpha+\frac \lambda 2]\subset
\left[\tau_1-N\cdot\frac{\lambda}2,
\tau_1+N\cdot\frac{\lambda}2\right]$$
 of the function
$\Psi^{(r-1)}(t)=\Psi^{(r-1)}_{0,a_2,b,\lambda}(t)$ the function
$\Delta_N^{(r-1)}(t)$ changes sign at least three times. But then the difference 
$$
\Psi^{(r-1)}_{0,0,b,\lambda}(t)-x^{(r-1)}_N(t)
$$
changes the sign at least three times on some monotonicity interval of the function $\Psi^{(r-1)}_{0,0,b,\lambda}(t)$ too. However this contradicts to the Kolmogorov comparison theorem (see theorem~A and, for example,~\cite[Statement 5.5.3]{korn2}) because the Euler spline $\Psi^{(r-1)}_{0,0,b,\lambda}(t)$ is comparison function for the function $x^{(r-1)}_N(t)$.

If the condition $b)$ of the theorem holds, then applying similar arguments we will get contradiction with Kolmogorov comparison theorem.

If the condition $c)$ of the theorem holds, then applying similar arguments we will get contradiction with already proved case when condition $a)$ holds.
The theorem is proved.

\section{Some applications.}

From the theorem~\ref{th::1} we immediately get
\begin{lemma}\label{sl::1}
Let $r\in{\mathbb N}$, $x\in L_{\infty,\infty}^r({\mathbb R})$ and one of the conditions $a)-c)$ of the theorem~\ref{th::1} holds. Then on each monotonicity interval of the function $\Psi_{a_1,a_2,b,\lambda}(t)$ the difference $\Psi_{a_1,a_2,b,\lambda}(t)-x(t)$ has at most one sign change.
\end{lemma}
\noindent For 1-periodic non-negative integrable on period function $x(\cdot)$ denote by $r(x,\cdot)$ the decreasing rearrangement of the function $x$ (see, for example~\cite[Chapter~6]{korn1}).

As a corollary of the theorem~\ref{th::1} and the results of the Chapter~3 of the monograph~\cite{korn2} we get the following theorem.

\begin{theorem}\label{th::2}
Let $r\in{\mathbb N}$ and $1$--periodic function $x\in L^{r}_{\infty,\infty}({\mathbb R})$ are given. Let one of the conditions $a)-c)$ of the theorem~\ref{th::1} holds. Then for all $t>0$
$$\int\limits_0^tr(|x'|,u)du\leq \lambda^{r-1}\int\limits_0^tr(|\Psi'_{a_1,a_2,b,1}|,u)du.$$
\end{theorem}
For $a,b\in{\mathbb R}$, $a<b$, $p\in (0,\infty)$ and continuous function $x\colon {\mathbb R}\to{\mathbb R}$ set $\|x\|_{L_{p(a,b)}}:=\left(\int\limits_a^b\left|x(t)\right|^pdt\right)^{\frac{1}{p}}$.

From the theorem~\ref{th::2} and general theorems about rearrangements comparison (see, for example,~\cite[Statement~1.3.10]{korn3}) we get the following analogue of the Ligun inequality~\cite{Ligun} (see also~\cite{BKKP}, Chapter~6).

\begin{theorem}
Let $r\in{\mathbb N}$ and $1$--periodic function $x\in L^{r}_{\infty,\infty}({\mathbb R})$ are given. Let one of the conditions  $a)-c)$ of the theorem~\ref{th::1} holds. Then for all $1\leq p < \infty$ and natural  $k<r-2$ $($ if condition $a)$ holds, then for all natural $k<r-1)$
$$\|x^{(k)}\|_{L_{p(0,1)}}\leq \lambda^{r-k}\|\Psi_{a_1,a_2,b,1}^{(k)}\|_{L_{p(0,1)}}.$$
\end{theorem}

The following lemma is an analogue to the Bohr-Favard inequality, see, for example,~\cite{korn4}, Chapter~6.
\begin{lemma}\label{th::3}
Let $r\in{\mathbb N}$ and $1$--periodic function $x\in L^{r}_{\infty,\infty}({\mathbb R})$ are given. Let for $\lambda = 1$ one of the following conditions holds.
\begin{enumerate}
\item[$a)$] Numbers $a_1 = 0$, $a_2\geq 0$ and $b\neq 0$ are such that
 \begin{equation*}
\left\|x^{(k)}\right\|\le
\left\|\Psi^{(k)}_{a_1,a_2,b,\lambda}\right\|,\,\,k\in\{r-1,r\}.
 \end{equation*}
\item[$b)$] Numbers $a_1 \geq 0$, $a_2= 0$ and $b\neq 0$ are such that
 \begin{equation*}
\left\|x^{(k)}\right\|\le
\left\|\Psi^{(k)}_{a_1,a_2,b,\lambda}\right\|,\,\,k\in\{r-2,r\}.
 \end{equation*}
\item[$c)$] Numbers $a_1 \geq 0$, $a_2\geq 0$ and $b\neq 0$ are such that
 \begin{equation*}
\left\|x^{(k)}\right\|\le
\left\|\Psi^{(k)}_{a_1,a_2,b,\lambda}\right\|,\,\,k\in\{r-1,r-2,r\}.
 \end{equation*}
\end{enumerate}
Then $$E_0(x):=\inf\limits_{c\in{\mathbb R}}\|x-c\|\leq\left\|\Psi_{a_1,a_2,b,1}\right\|.$$
\end{lemma}

We will proceed by induction on $r$. The basis of the induction easily follows. We will dwell on the induction step.
Assume the contrary, let $E_0(x)>\left\|\Psi_{a_1,a_2,b,1}\right\|.$ Let $c$ be a constant of the best uniform approximation of the function $x$. We can count that $\max\limits_{t\in [0,1]}[x(t)-c]$ is attained in the point $t=0$, $\min\limits_{t\in [0,1]}[x(t)-c]$ --- in the point $m$ and
\begin{equation}\label{000}
m<\frac 12.
\end{equation}
This means that $x'(0) = x'(m)=0$. Moreover $$-\int\limits_0^mx'(t)dt=x(0)-x(m)=2E_0(x)>2\left\|\Psi_{a_1,a_2,b,1}\right\|=-\int\limits_0^m\Psi_{a_1,a_2,b,1}'(t)dt.$$
However the last inequality together with the induction hypothesis and~\eqref{000} contradicts the lemma~\ref{sl::1}.

Using the theorem~\ref{th::1}, lemma~\ref{th::3} and ideas from~\cite{BKP1} (see also \S~6.4 of the monograph~\cite{BKKP}) we get the following analogue of Babenko, Kofanov and Pichugov inequality.
\begin{theorem}\label{th::4}

Let $r\in{\mathbb N}$ and $1$--periodic function $x\in L^{r}_{\infty,\infty}({\mathbb R})$ are given. Let for some $\lambda > 0$ one of the conditions $a)-c)$ of the lemma~\ref{th::3} holds
 and $E_0(x)=\left\|\Psi_{a_1,a_2,b,\lambda}\right\|.$
Then $$\|x\|_{L_{p(0,1)}}\geq\|\Psi_{a_1,a_2,b,\lambda}\|_{L_{p(0,\lambda)}}.$$
\end{theorem}

For a function $x\in L_\infty({\mathbb R})$ let $c(x)$ denote the constant of the best approximation for the function $x$ in $L_\infty({\mathbb R})$.

From theorem~\ref{th::4}, theorem~\ref{th::1} and ideas from~\cite{BKP2} (see also \S~6.7 of the monograph~\cite{BKKP}) we get the following analogue of Nagy type inequality (see~\cite{Nad1}), that was obtained by Babenko, Kofanov and Pichugov.
\begin{theorem}
Let $r\in{\mathbb N}$, $1$--periodic function $x\in L^{r}_{\infty,\infty}({\mathbb R})$, and numbers $p,q\in (0,\infty)$, $q>p$ be given.
Let for some $\lambda > 0$ one of the conditions $a)-c)$ of the lemma~\ref{th::3} hold and
$\left\|\Psi_{a_1,a_2,b,\lambda}\right\|_{L_{p(0,\lambda)}}=\|x-c(x)\|_{L_{p(0,1)}}.$ Then $$\left\|\Psi_{a_1,a_2,b,\lambda}\right\|_{L_{q(0,\lambda)}}\geq\|x-c(x)\|_{L_{q(0,1)}}.$$
\end{theorem}

\end{document}


{\bf \Large 517.5}
\begin{center}
{\Large \bf Теоремы сравнения производных и некоторые их приложения}
\vskip 5 mm.
{\Large\bf  В.~Ф.~Бабенко, О.~В.~Коваленко}
\vskip 5 mm.
{Днепропетровский национальный университет
им. Олеся Гончара,\\
Днепропетровск 49050. \\{\it E-mail: babenko.vladislav@gmail.com }
{\it E-mail: olegkovalenko90@gmail.com} }
\end{center}

\begin{abstract}
{Получены аналоги теорем сравнения Колмогорова и указаны их некоторые  приложения.}
\end{abstract}

\section{Обозначения. Постановка задачи. Известные результаты.} Через
$L_\infty({\mathbb R})$ будем обозначать пространство измеримых и
сущест\-вен\-но ограничен\-ных функций $x\colon {\mathbb R}\to {\mathbb R}$ с нормой
$$\|x\|=\|x\|_{L_\infty({\mathbb R})}={\rm
ess\,sup}\left\{|x(t)|:t\in{\mathbb R}\right\}.$$ Для натурального $r$ через
$L_\infty^r({\mathbb R})$ обозначим пространство  функций $x\colon {\mathbb R}\to
{\mathbb R}$ таких, что производная $x^{(r-1)}, \; x^{(0)}=x,$ локально
абсолютно непрерывна, и $x^{(r)}\in L_\infty({\mathbb R})$. Пусть также
$L_{\infty,\infty}^r({\mathbb R}):=L_\infty^r({\mathbb R})\bigcap L_\infty({\mathbb R})$.

Для $r\in {\mathbb N}$ через $\varphi_r(t)$ будем обозначать сплайн Эйлера порядка $r$ (т.е. $r$-ю периодическую первообразную функции ${\rm sgn}\sin t$ со средним значением нуль на периоде).
Для $\lambda> 0$ положим $\varphi_{\lambda, r}(t):=\lambda^{-r}\varphi_{r}(\lambda t).$

Для доказательства своего знаменитого неравенства (см. \cite{Kol1}---\cite{Kol3}) Колмогоровым было доказано
утверж\-дение, известное как теорема сравнения.

{\bf Теорема A. }{\it Пусть $r\in{\mathbb N}$ и задана функция $x\in L^r_{\infty,\infty}({\mathbb R})$. Пусть числа $a\in{\mathbb R}$ и $\lambda>0$ таковы, что
$$
\|x^{(k)}\|\leq \|a\varphi_{\lambda,r}^{(k)}\|,\,\, k\in\{0,r\}.
$$ Если точки $\xi,\eta\in{\mathbb R}$ таковы, что $x(\xi)=a\varphi_{\lambda,r}(\eta)$, то $$|x'(\xi)|\leq |a|\cdot|\varphi'_{\lambda,r}(\eta)|.$$ }
Как сама теорема сравнения Колмогорова, так и метод её
доказательства сыграли большую роль при точном решении многих
экстремальных задач теории приближений (см.~\cite{korn1, korn2, korn3, BKKP}).

Цель данной статьи -- доказательство некоторых аналогов теоремы сравнения Колмогорова.

В следующем параграфе мы введем семейство сплайнов, которые в дальнейшем будут играть роль, которую в теореме А играли сплайны Эйлера, и
изучим некоторые их свойства.  В параграфе 3 мы докажем три аналога теоремы
сравнения Колмогорова для случаев, когда заданы норма функции и нормы
ее произ\-вод\-ных порядков $r-1$ и $r$; норма функции и нормы её производных порядков $r-2$ и $r$; норма функции и нормы
ее произ\-вод\-ных порядков $r-2$, $r-1$ и $r$. В параграфе 4 мы приведем некоторые приложения полученных теорем сравнения.

\medskip

\section{Функции сравнения и их свойства.} Пусть $a_1, a_2\geq 0$.  Положим $T:=a_1+a_2+2$. Определим
функ\-цию $\psi_1(a_1,a_2;t)$ следующим образом. На отрезке $[0,T]$
положим
 $$
           \psi_1(a_1,a_2;t):=\left\{
           \begin{array}{rcl}
            0, &  t\in[0,a_1], \\
              t-a_1,   & t\in[a_1,a_1+1], \\
              1, & t\in [a_1+1,a_1+a_2+1], \\
              2+a_1+a_2-t, & t\in [a_1+a_2+1, T].\\
           \end{array}
           \right.
  $$
  Продолжим функцию $\psi_1(a_1,a_2;t)$ на отрезок $[T,2\cdot T]$ равенством
  \begin{equation}\label{2.0}
  \psi_1(a_1,a_2;t) = - \psi_1(a_1,a_2;t-T),\,t\in[T,2\cdot T]
  \end{equation} а затем
  периодически с периодом $2\cdot T$ на всю ось.

  Заметим, что $\psi_1(a_1,a_2;t)\in L_{\infty,\infty}^1({\mathbb R})$.
  Для $r\in {\mathbb N}$ обозначим через $\psi_r(a_1,a_2;t)$ $(r-1)$-ю $2\cdot T$ -- периодическую
  первообразную функции $\psi_1(a_1,a_2;t)$ с нулевым средним на периоде (так что, в частности,
  $\psi'_r(a_1,a_2;t)=\psi_{r-1}(a_1,a_2;t)$).
Функции $\psi_r(a_1,a_2;t)$ впервые рассматривал Родов~\cite{Rodov1946}.

  Приведем некоторые свойства функции $\psi_r(a_1,a_2;t)$, $r\in{\mathbb N}$, которые
  нетруд\-но уста\-но\-вить либо непосредственно, либо по аналогии со
  свойствами эйлеровых иде\-альных сплайнов $\varphi_r$ (см., например,~\cite[гл. 5]{korn1},~\cite[гл. 3]
  {korn2}). Отметим, что функция $\psi_r$ имеет период $2\cdot T$ и для всех $r\geq 1$ 
   $$\psi_r(a_1,a_2;t) = - \psi_r(a_1,a_2;t-T),\,t\in[T,2\cdot T].$$
    Кроме того функция $\psi_2(a_1,a_2;t)$ имеет ровно два нуля на периоде -- точки $a_1+\frac{a_2}{2}+1$ и $2a_1+\frac{3a_2}{2}+3$.
    Следовательно функция   $\psi_r(a_1,a_2;t)$ при $r\geq 2$ также имеет ровно два нуля на
    периоде: для любого $k\in {\mathbb N}$
    \begin{equation}\label{1}
    \psi_{2k+1}(a_1,a_2;0)=\psi_{2k+1}(a_1,a_2;a_1+a_2+2)=0,
    \end{equation}
        \begin{equation}\label{2}
     \psi_{2k}\left(a_1,a_2;a_1+\frac{a_2}{2}+1\right)=\psi_{2k}\left(a_1,a_2;2a_1+\frac{3a_2}{2}+3\right)=0.
    \end{equation}

Отметим, что при $a_1 = 0$ равенство~\eqref{1} справедливо и для $k=0$.

Отсюда, в свою очередь, следует, что при $r\geq 3$ (а в случае, когда $a_1=0$ для $r\geq 2$) функция
$\psi_r(a_1,a_2;t)$ строго монотонна между нулями своей
 производной,
 а график функции $\psi_r(a_1,a_2;t)$ является выпуклым на каждом промежутке
 знакопостоянства. Кроме того, как легко видеть, график
 $\psi_r(a_1,a_2;t)$ симметричен относительно ее нулей, а также
 относительно прямых вида $t=t_0$, где $t_0$ -- нуль $\psi'_r(a_1,a_2;t)$.
 Наконец отметим, что $\psi_r(0,0;t)=\varphi_{\pi/2, r}(t)$.

Для $r\in{\mathbb N}$, $a_1,a_2\geq 0,\,\lambda>0$ и $b\in{\mathbb R}$ положим $$\Psi_{a_1,a_2,b,\lambda}(t)=\Psi_{r;a_1,a_2,b,\lambda}(t):=$$
$$=b\left(\frac{\lambda}{2a_1+2a_2+4}\right)^{r}\psi_r\left(a_1,a_2;\frac{2a_1+2a_2+4}{\lambda} t\right).$$ Отметим, что функция $\Psi_{a_1,a_2,b,\lambda}(t)$ является $\lambda$ -- периодической.

 \begin{theorem}\label{th::0}
 Пусть $r\in{\mathbb N}$ и $x\in L^r_{\infty,\infty}({\mathbb R})$. Тогда
 \begin{enumerate}
 \item[$a)$]   существуют $a_2\geq 0,\,\lambda>0$ и $b\in{\mathbb R}$ такие, что
 $$ \left\|\Psi_{0,a_2,b,\lambda}^{(s)}\right\|=\left\|x^{(s)}\right\|,\,\,s\in\left\{0,r-1,r\right\}.$$
  \item[$b)$]  существуют $a_1\geq 0,\,\lambda>0$ и $b\in{\mathbb R}$ такие, что
  $$ \left\|\Psi_{a_1,0,b,\lambda}^{(s)}\right\|=\left\|x^{(s)}\right\|,\,\,s\in\left\{0,r-2,r\right\}.$$
   \item[$c)$] существуют $a_1,a_2\geq 0,\,\lambda>0$ и $b\in{\mathbb R}$ такие, что
    $$ \left\|\Psi_{a_1,a_2,b,\lambda}^{(s)}\right\|=\left\|x^{(s)}\right\|,\,\,s\in\left\{0,r-1,r-2,r\right\}.$$
 \end{enumerate}
\end{theorem}


Справедливость теоремы~\ref{th::0} по сути следует из теоремы сравнения Колмогорова. Докажем утверждение $a)$, остальные утверждения доказываются аналогично.
 
 Ясно, что $\Psi_{0,0,b,\lambda}(t)=b\left(\frac\lambda 4\right)^r\varphi_{\pi/2}(\frac 4\lambda t)$. Следовательно, параметры $b$ и $\lambda$ можно выбрать
   так, чтобы $ \left\|\Psi_{0,0,b,\lambda}^{(s)}\right\|=\left\|x^{(s)}\right\|$, $s=r-1,r$. Тогда из теоремы A будет следовать, что $ \left\|\Psi_{0,0,b,\lambda}\right\|\leq\left\|x\right\|$. При возрастании параметра $a_2$ $ \left\|\Psi_{0,a_2,b,\lambda}\right\|$ будет непрерывно изменяясь возрастать от $ \left\|\Psi_{0,0,b,\lambda}\right\|$ до $\infty$, а $ \left\|\Psi_{0,a_2,b,\lambda}^{(s)}\right\|,\; s=r-1, r$ будут оставаться неизменными. Поэтому можно выбрать параметр $a_2$ так, чтобы $ \left\|\Psi_{0,a_2,b,\lambda}^{(s)}\right\|=\left\|x^{(s)}\right\|$, $s=0,r-1,r$.


\section{Теоремы сравнения.}

Следу\-ющая теорема содержит в себе три аналога теоремы сравнения Колмогорова.
 \begin{theorem}\label{th::1}
 Пусть $r\in{\mathbb N}$ и $x\in L_{\infty,\infty}^r({\mathbb R})$. Пусть выполняется одно из следующих условий.
\begin{enumerate}
\item[$a)$] Числа $a_1 = 0$, $a_2\geq 0$, $\lambda>0$ и $b\neq 0$ таковы, что
 \begin{equation}\label{4.1}
\left\|x^{(k)}\right\|\le
\left\|\Psi^{(k)}_{a_1,a_2,b,\lambda}\right\|,\,\,k\in\{0,r-1,r\}.
 \end{equation}
\item[$b)$] Числа $a_1 \geq 0$, $a_2= 0$, $\lambda>0$ и $b\neq 0$ таковы, что
 \begin{equation}\label{4.2}
\left\|x^{(k)}\right\|\le
\left\|\Psi^{(k)}_{a_1,a_2,b,\lambda}\right\|,\,\,k\in\{0,r-2,r\}.
 \end{equation}
\item[$c)$] Числа $a_1 \geq 0$, $a_2\geq 0$, $\lambda>0$ и $b\neq 0$ таковы, что
 \begin{equation}\label{4.3}
\left\|x^{(k)}\right\|\le
\left\|\Psi^{(k)}_{a_1,a_2,b,\lambda}\right\|,\,\,k\in\{0,r-1,r-2,r\}.
 \end{equation}
\end{enumerate}

  Если точки $\tau$ и $\xi$ таковы, что $x(\tau) = \Psi_{a_1,a_2,b,\lambda}(\xi)$, то
\begin{equation}\label{6'}\left|x'(\tau)\right|\le\left|\Psi'_{a_1,a_2,b,\lambda}(\xi)\right|.
 \end{equation}
 \end{theorem}

{\bf Доказательство.} Для сокращения записей в ходе данного
доказательства мы будем писать $\Psi (t)$ вместо
$\Psi_{a_1,a_2,b,\lambda}(t)$. Рассматривая при необходимости
функ\-цию $-x(t)$ вместо функции $x(t)$ и функцию $-\Psi (t)$ вместо
$\Psi (t)$, мы можем считать, что $x'(\tau)>0$ и
 \begin{equation}\label{7}
 \Psi'(\tau)>0.
 \end{equation}
Кроме того, рассматривая подходящий сдвиг $\Psi (\cdot +\alpha)$
функции $\Psi$ можем считать, что  и $\tau = \xi$, т.е.
 \begin{equation}\label{6}
 x(\tau)=\Psi(\tau).
 \end{equation}
Предположим, что~(\ref{6}) имеет место, но при этом вместо
неравенства~(\ref{6'}) (с $\xi =\tau$) имеет место неравенство
$$\left|x'(\tau)\right|>\left|\Psi'(\tau)\right|. $$
 Обозначим через $(\tau_1,\tau_2)$ наименьший интервал монотонности $\Psi$, содер\-жа\-щий
 точку $\tau$ и такой, что $\Psi'(\tau_1)=\Psi'(\tau_2)=0$.
В силу сделанного предположе\-ния существует число $\delta>0$ такое,
что $x'(t)>\Psi'(t)$ для всех $t\in (\tau-\delta,\tau+\delta)$, а
значит в силу \eqref{6} $x(\tau+\delta)>\Psi(\tau+\delta)$ и
$x(\tau-\delta)<\Psi(\tau-\delta)$.

Выберем $\varepsilon>0$ настолько малым, чтобы для функции
$x_\varepsilon(t):=(1-\varepsilon)x(t)$ выполнялись неравенства
$x_\varepsilon(\tau+\delta)>\Psi(\tau+\delta)$ и
$x_\varepsilon(\tau-\delta)<\Psi(\tau-\delta)$. В силу одного из соотношений~\eqref{4.1} -- \eqref{4.3} и
\eqref{7} будет
$$x_\varepsilon(\tau_1)>\Psi(\tau_1),\;\; x_\varepsilon(\tau_2)<\Psi(\tau_2).$$
Таким образом на промежутке $(\tau_1,\tau_2)$ разность
$\Delta_\varepsilon(t):=x_\varepsilon(t)-\Psi(t)$ будет иметь не
менее трёх перемен знака.

Как легко видеть,  существует последовательность функций $\mu_N\in
C^\infty({\mathbb R})$, $N\in {\mathbb N}$, со следующими свойствами:
\begin{enumerate}
\item $\mu_N(t)=1$ на промежутке $[\tau_1,\tau_2]$; $\|\mu_N\|=1$;
\item $\mu_N(t)=0$ для всех $t$ вне промежутка
$[\tau_1-N\cdot\frac{\lambda}2;\tau_1+N\cdot\frac{\lambda}2]$ (где, как и ранее, $T=a_1+a_2+2$);
\item для всех
$k=1,2,\dots,r$
$$\max\limits_{j=\overline{1,k}}\left\|\mu_N^{(j)}\right\|<
\varepsilon \|x_\varepsilon^{(k)}\|\left(\sum\limits_{i=1}^k
C_k^i\left\|x_\varepsilon^{(k-i)}\right\|\right)^{-1},$$ если $N$
достаточно велико.
\end{enumerate}

Ниже считаем, что $N$ выбрано настолько большим, что свойство 3
выполнено.

Положим $$x_N(t):=x_\varepsilon(t)\cdot\mu_N(t),$$
и
$$\Delta_N(t):=\Psi(t)-x_N(t).$$ Тогда
$$x_N(t)=x_\varepsilon(t),\;\textcyrillic{если }\; t\in[\tau_1,\tau_2],$$
\begin{equation}\label{7'}
\Delta_N(t)=\Psi(t),\textcyrillic{если }\,|t-\tau_1|\geq N\cdot\frac\lambda 2
\end{equation}
и
$$
\|x_N\|\leq\|x_\varepsilon\|=(1-\varepsilon)\|x\|\leq(1-\varepsilon)\|\Psi\|.
$$
Кроме того для  $k=1,\dots,r$ $$\left|x_N^{(k)}(t)\right|=
\left|\left[x_\varepsilon(t)\mu_N(t)\right]^{(k)}\right|=
\left|\sum\limits_{i=0}^k C_k^i
x_\varepsilon^{(k-i)}(t)\mu_N^{(i)}(t)\right|\leq
$$
$$
\leq \left\|x_\varepsilon^{(k)}\right\|+\sum\limits_{i=1}^k
C_k^i\left\|x_\varepsilon^{(k-i)}\right\|\left\|\mu_N^{(i)}\right\|.$$ Отсюда,
с учетом свойства 3 функции $\mu_N$ и выбора числа $N$, получаем
$$
\left\|x_N^{(k)}\right\|<
\left\|x_\varepsilon^{(k)}\right\|+\varepsilon\left\|x^{(k)}_\varepsilon\right\|=(1-\varepsilon)\left\|x^{(k)}\right\|+\varepsilon\left\|x^{(k)}\right\|=\left\|x^{(k)}\right\|.
$$
Для $t\in [\tau_1,\tau_2]$ имеем
$\Delta_N=\Psi(t)-x_\varepsilon(t),$ а значит функция $\Delta_N(t)$
имеет не менее  трех перемен знака на отрезке $[\tau_1,\tau_2]$. На
каждом из остальных промежутков монотонности функции $\Psi$ функция
$\Delta_N$ имеет не менее одной перемены знака. Таким образом на
промежутке
$\left[\tau_1-N\cdot\frac{\lambda}2,\tau_1+N\cdot\frac{\lambda}2\right]$
функция $\Delta_N(t)$ имеет не менее $2N+2$ перемен знака. Кроме
того, в силу \eqref{1},\eqref{2} и \eqref{7'} для всех
$i=1,2,\dots,\left[\frac{r-1}{2}\right]$ справедливы равенства
\begin{equation}\label{10}
\Delta_N^{(2i-1)}\left(\tau_1-N\cdot\frac\lambda 2\right)=\Delta_N^{(2i-1)}
\left(\tau_1+N\cdot\frac\lambda 2\right)=0.
\end{equation}

Все проведённые выше рассуждения справедливы при выполнении любого из условий $a)-c)$.
Пусть теперь выполняется условие $a)$ теоремы.

Применяя теорему Ролля и учитывая \eqref{10} получаем, что
функция $\Delta_N^{(r-1)}(t)$ имеет не менее $2N+2$ нуля на
промежутке
$$\left[\tau_1-N\cdot\frac{\lambda} 2,\tau_1+N\cdot\frac{\lambda} 2\right].$$
Отсюда следует, что на некотором промежутке  монотонности
$$[\alpha,\alpha+\frac \lambda 2]\subset
\left[\tau_1-N\cdot\frac{\lambda}2,
\tau_1+N\cdot\frac{\lambda}2\right]$$
 функции
$\Psi^{(r-1)}(t)=\Psi^{(r-1)}_{0,a_2,b,\lambda}(t)$ функция
$\Delta_N^{(r-1)}(t)$ меняет знак не менее трёх раз. Но тогда и
разность
$$
\Psi^{(r-1)}_{0,0,b,\lambda}(t)-x^{(r-1)}_N(t)
$$
на некотором промежутке монотонности функции $\Psi^{(r-1)}_{0,0,b,\lambda}(t)$ меняет знак не менее трех раз. Однако это противоречит
теореме сравнения Колмогорова (см., теорему~A и, например,~\cite[Предложение
5.5.3]{korn2}) так как эйлеров сплайн $\Psi^{(r-1)}_{0,0,b,\lambda}(t)$ является функцией сравнения для функции $x^{(r-1)}_N(t)$.

Если выполняется условие $b)$ теоремы, то применяя аналогичные рассуждение придём к противоречию с теоремой сравнения Колмогорова.

Если же выполняется условие $c)$ теоремы, то применяя аналогичные рассуждения придём к противоречию с уже доказанным случаем $a)$.
Теорема доказана.

\section{Некоторые приложения.}

Из теоремы~\ref{th::1} сразу получаем, что справедлива
\begin{lemma}\label{sl::1}
Пусть $r\in{\mathbb N}$, $x\in L_{\infty,\infty}^r({\mathbb R})$ и выполняется одно из  условий $a)-c)$ теоремы~\ref{th::1}. Тогда на каждом промежутке монотонности функции $\Psi_{a_1,a_2,b,\lambda}(t)$ разность $\Psi_{a_1,a_2,b,\lambda}(t)-x(t)$ имеет не более одной перемены знака.
\end{lemma}
\noindent Для 1-периодической неотрицательной суммируемой на периоде функции $x(\cdot)$ через $r(x,\cdot)$ будем обозначать убывающую перестановку функции $x$ (см., например~\cite{korn1}, гл. 6).

Как следствие из теоремы~\ref{th::1} и результатов главы~3 монографии~\cite{korn2} получаем следующую теорему.

\begin{theorem}\label{th::2}
Пусть заданы $r\in{\mathbb N}$ и $1$--периодическая функция $x\in L^{r}_{\infty,\infty}({\mathbb R})$. Пусть выполнено одно из условий $a)-c)$ теоремы~\ref{th::1}. Тогда для всех $t>0$
$$\int\limits_0^tr(|x'|,u)du\leq \lambda^{r-1}\int\limits_0^tr(|\Psi'_{a_1,a_2,b,1}|,u)du.$$
\end{theorem}
Для $a,b\in{\mathbb R}$, $a<b$, $p\in (0,\infty)$ и  непрерывной функции $x\colon {\mathbb R}\to{\mathbb R}$ положим $\|x\|_{L_{p(a,b)}}:=\left(\int\limits_a^b\left|x(t)\right|^pdt\right)^{\frac{1}{p}}$.

Из теоремы~\ref{th::2} и общих теорем о сравнении перестановок (см., например,~\cite{korn3}, \textcyrillic{предложение} 1.3.10) получаем следующий аналог неравенства Лигуна~\cite{Ligun} (см. также~\cite{BKKP}, глава 6).

\begin{theorem}
Пусть заданы $r\in{\mathbb N}$ и $1$--периодическая функция $x\in L^{r}_{\infty,\infty}({\mathbb R})$. Пусть выполнено одно из условий $a)-c)$ теоремы~\ref{th::1}. Тогда для любого $1\leq p < \infty$ и натурального $k<r-2$ $($а если выполняется условие $a)$, то для любого натурального $k<r-1)$
$$\|x^{(k)}\|_{L_{p(0,1)}}\leq \lambda^{r-k}\|\Psi_{a_1,a_2,b,1}^{(k)}\|_{L_{p(0,1)}}.$$
\end{theorem}

Следующая лемма является аналагом неравенства Бора-Фавара, см., напр.,~\cite{korn4}, гл. 6.
\begin{lemma}\label{th::3}
Пусть заданы $r\in{\mathbb N}$ и $1$--периодическая функция $x\in L^{r}_{\infty,\infty}({\mathbb R})$. Пусть для $\lambda = 1$ выполняется одно из следующих условий.
\begin{enumerate}
\item[$a)$] Числа $a_1 = 0$, $a_2\geq 0$ и $b\neq 0$ таковы, что
 \begin{equation*}
\left\|x^{(k)}\right\|\le
\left\|\Psi^{(k)}_{a_1,a_2,b,\lambda}\right\|,\,\,k\in\{r-1,r\}.
 \end{equation*}
\item[$b)$] Числа $a_1 \geq 0$, $a_2= 0$ и $b\neq 0$ таковы, что
 \begin{equation*}
\left\|x^{(k)}\right\|\le
\left\|\Psi^{(k)}_{a_1,a_2,b,\lambda}\right\|,\,\,k\in\{r-2,r\}.
 \end{equation*}
\item[$c)$] Числа $a_1 \geq 0$, $a_2\geq 0$ и $b\neq 0$ таковы, что
 \begin{equation*}
\left\|x^{(k)}\right\|\le
\left\|\Psi^{(k)}_{a_1,a_2,b,\lambda}\right\|,\,\,k\in\{r-1,r-2,r\}.
 \end{equation*}
\end{enumerate}
Тогда $$E_0(x):=\inf\limits_{c\in{\mathbb R}}\|x-c\|\leq\left\|\Psi_{a_1,a_2,b,1}\right\|.$$
\end{lemma}

Доказательство проведем индукцией по $r$. Базис индукции проверяется непосредственно. Остановимся на индуктивном шаге.
Предположим противное, пусть $E_0(x)>\left\|\Psi_{a_1,a_2,b,1}\right\|.$ Пусть $c$ --- константа наилучшего равномерного приближения функции $x$. Можем считать, что $\max\limits_{t\in [0,1]}[x(t)-c]$ достигается в точке $t=0$, $\min\limits_{t\in [0,1]}[x(t)-c]$ --- в точке $m$ и при этом
\begin{equation}\label{000}
m<\frac 12.
\end{equation}
Это значит, что $x'(0) = x'(m)=0$. Кроме того $$-\int\limits_0^mx'(t)dt=x(0)-x(m)=2E_0(x)>2\left\|\Psi_{a_1,a_2,b,1}\right\|=-\int\limits_0^m\Psi_{a_1,a_2,b,1}'(t)dt.$$
Однако последнее неравенство вместе с индуктивным предположением и~\eqref{000} противоречит лемме~\ref{sl::1}.

Используя теорему~\ref{th::1}, лемму~\ref{th::3} и идеи из~\cite{BKP1} (см. также параграф 6.4 монографии~\cite{BKKP}) получаем следующий аналог неравенства Бабенко, Кофанова и Пичугова.
\begin{theorem}\label{th::4}

Пусть заданы $r\in{\mathbb N}$ и $1$--периодическая функция $x\in L^{r}_{\infty,\infty}({\mathbb R})$. Пусть для некоторого $\lambda > 0$ выполняется одно из условий $a)-c)$ леммы~\ref{th::3}
 и $E_0(x)=\left\|\Psi_{a_1,a_2,b,\lambda}\right\|.$
Тогда $$\|x\|_{L_{p(0,1)}}\geq\|\Psi_{a_1,a_2,b,\lambda}\|_{L_{p(0,\lambda)}}.$$
\end{theorem}

Для функции $x\in L_\infty({\mathbb R})$ через $c(x)$ обозначим константу наилучшего приближения функции $x$ в $L_\infty({\mathbb R})$.

Из теоремы~\ref{th::4}, теоремы~\ref{th::1} и идей из~\cite{BKP2} (см. также параграф 6.7 монографии~\cite{BKKP}) получаем следующий аналог неравенства типа Надя (см.~\cite{Nad1}), полученного Бабенко, Кофановым и Пичуговым.
\begin{theorem}
Пусть заданы $r\in{\mathbb N}$, $1$--периодическая функция $x\in L^{r}_{\infty,\infty}({\mathbb R})$, и числа $p,q\in (0,\infty)$, $q>p$.
Пусть для некоторого $\lambda > 0$ выполняется одно из условий $a)-c)$ леммы~\ref{th::3} и
$\left\|\Psi_{a_1,a_2,b,\lambda}\right\|_{L_{p(0,\lambda)}}=\|x-c(x)\|_{L_{p(0,1)}}.$ Тогда $$\left\|\Psi_{a_1,a_2,b,\lambda}\right\|_{L_{q(0,\lambda)}}\geq\|x-c(x)\|_{L_{q(0,1)}}.$$
\end{theorem}